\documentclass{amsproc}

\usepackage{}

\newtheorem{theorem}{Theorem}[section]

\newtheorem{conjecture}[theorem]{Conjecture}
\newtheorem{corollary}[theorem]{Corollary}

\newtheorem{lemma}[theorem]{Lemma}

\newtheorem{problem}[theorem]{Problem}

\theoremstyle{definition}
\newtheorem{definition}[theorem]{Definition}

\theoremstyle{remark}

\numberwithin{equation}{section}

\begin{document}

\title{PLANK THEOREMS VIA SUCCESSIVE INRADII}


\author{K\'aroly Bezdek}
\address{Department of Mathematics and\\Statistics\\
University of Calgary\\
Calgary, Alberta, Canada T2N 1N4}
\email{bezdek@math.ucalgary.ca}
\thanks{Partially supported by a Natural Sciences and Engineering Research Council of Canada Discovery Grant.}


\subjclass[2010]{Primary 52C17, 05B40, 11H31, and 52C45.}

\date{June 27, 2013}

\begin{abstract}
In the 1930's, Tarski introduced his plank problem at a time when the field discrete geometry was about to born. It is quite remarkable that
Tarski's question and its variants continue to generate interest in the geometric as well as analytic aspects of coverings by planks in the present
time as well. Besides giving a short survey on the status of the affine plank conjecture of Bang (1950) we prove some new partial results 
for the successive inradii of the convex bodies involved. The underlying geometric structures are successive hyperplane cuts introduced several 
years ago by Conway and inductive tilings introduced recently by Akopyan and Karasev.
\end{abstract}

\maketitle


\section{Introduction}

As usual, a {\it convex body} of the Euclidean space $\mathbb{E}^d$ is a compact convex set with non-empty interior. Let $\mathbf{C}\subset\mathbb{E}^d$ be a convex body, and let $H\subset\mathbb{E}^d$ be a hyperplane. Then the distance $w(\mathbf{C} , H)$ between the two supporting hyperplanes of $\mathbf{C}$ parallel to $H$ is called the {\it width of $\mathbf{C}$ parallel to $H$}. Moreover, the smallest width of $\mathbf{C}$ parallel to hyperplanes of $\mathbb{E}^d$ is called the {\it minimal width} of $\mathbf{C}$ and is denoted by $w(\mathbf{C})$.

Recall that in the 1930's, Tarski posed what came to be known as the plank problem. A {\it plank} $\mathbf{P}$ in $\mathbb{E}^d$ is the (closed) set of points between two distinct parallel hyperplanes. The {\it width} $w(\mathbf{P})$ of $\mathbf{P}$ is simply the distance between the two boundary hyperplanes of $\mathbf{P}$. Tarski conjectured that if a convex body of minimal width $w$ is covered by a collection of planks in $\mathbb{E}^d$, then the sum of the widths of these planks is at least $w$. This conjecture was proved by Bang in his memorable paper \cite{Ba51}. (In fact, the proof presented in that paper is a simplification and generalization of the proof published by Bang somewhat earlier in \cite{Ba50}.) Thus, we call the following statement Bang's plank theorem.

\begin{theorem}\label{Bang-plank-th}
If the convex body $\mathbf{C}$ is covered by the planks $\mathbf{P}_1, \mathbf{P}_2, \dots , \mathbf{P}_n$ in $\mathbb{E}^d, d\ge 2$ (i.e., $\mathbf{C}\subset \mathbf{P}_1\cup \mathbf{P}_2\cup \dots \cup \mathbf{P}_n\subset\mathbb{E}^d$), then
$\sum_{i=1}^n w(\mathbf{P}_i)\ge w(\mathbf{C})$.
\end{theorem}

In \cite{Ba51}, Bang raised the following stronger version of Tarski's plank problem called the affine plank problem. We phrase it via the following definition.
Let $\mathbf{C}$ be a convex body and let $\mathbf{P}$ be a plank with boundary hyperplanes parallel to the hyperplane $H$ in $\mathbb{E}^d$. We define the {\it $\mathbf{C}$-width} of the plank $\mathbf{P}$ as $\frac{w(\mathbf{P}) }{w(\mathbf{C} , H) }$ and label it $w_{\mathbf{C}}(\mathbf{P})$. (This notion was introduced by Bang \cite{Ba51} under the name ``relative width''.)

\begin{conjecture}\label{Bang-conjecture}
If the convex body $\mathbf{C}$ is covered by the planks $\mathbf{P}_1, \mathbf{P}_2, \dots ,$ $\mathbf{P}_n$ in $\mathbb{E}^d, d\ge 2$, then
$\sum_{i=1}^n w_{\mathbf{C}}(\mathbf{P}_i)\ge 1$.
\end{conjecture}

The special case of Conjecture \ref{Bang-conjecture}, when the convex body to be covered is centrally symmetric, has been proved by Ball in \cite{Bal91}. Thus, the following is Ball's plank theorem.

\begin{theorem}\label{Ball-plank-th}
If the centrally symmetric convex body $\mathbf{C}$ is covered by the planks $\mathbf{P}_1, \mathbf{P}_2, \dots , \mathbf{P}_n$ in $\mathbb{E}^d, d\ge 2$, then
$\sum_{i=1}^n w_{\mathbf{C}}(\mathbf{P}_i)\ge 1$.
\end{theorem}

It was Alexander \cite{Al68} who noticed that Conjecture \ref{Bang-conjecture} is equivalent to the following
generalization of a problem of Davenport.

\begin{conjecture}\label{Alexander--Davenport}
If a convex body $\mathbf{C}$ in $\mathbb{E}^d, d\ge 2$ is sliced by $n-1$ hyperplane cuts, then there exists a piece that covers a translate of $\frac{1}{n}\mathbf{C}$.
\end{conjecture}

We note that the paper \cite{BeBe96} of A. Bezdek and the author proves Conjecture \ref{Alexander--Davenport} for successive hyperplane cuts (i.e., for hyperplane cuts when each cut divides one piece). Also, the same paper (\cite{BeBe96}) introduced two additional equivalent versions of Conjecture \ref{Bang-conjecture}. As they seem to be of independent interest we recall them following the terminology used in \cite{BeBe96}.

Let $\mathbf{C}$ and $\mathbf{K}$ be convex bodies in $\mathbb{E}^d$ and let $H$ be a hyperplane of $\mathbb{E}^d$. The {\it $\mathbf{C}$-width of $\mathbf{K}$ parallel to $H$} is denoted by $ w_{\mathbf{C}}(\mathbf{K} , H)$ and is defined as $\frac{w(\mathbf{K} , H)}{w(\mathbf{C} , H) }$. The {\it minimal $\mathbf{C}$-width of $\mathbf{K}$} is denoted by $ w_{\mathbf{C}}(\mathbf{K})$ and is defined as the minimum of $ w_{\mathbf{C}}(\mathbf{K} , H)$, where the minimum is taken over all possible hyperplanes $H$ of $\mathbb{E}^d$. Recall that the inradius of $\mathbf{K}$ is the radius of the largest ball contained in $\mathbf{K}$. It is quite natural then to introduce the {\it $\mathbf{C}$-inradius of $\mathbf{K}$} as the factor of the largest positive homothetic copy of $\mathbf{C}$, a translate of which is contained in $\mathbf{K}$. We need to do one more step to introduce the so-called successive $\mathbf{C}$-inradii of $\mathbf{K}$ as follows. 

Let $r$ be the $\mathbf{C}$-inradius of $\mathbf{K}$. For any $0<\rho\le r$ let the {\it $\rho\mathbf{C}$-rounded body of $\mathbf{K}$} be denoted by ${\mathbf{K}}^{\rho\mathbf{C}}$ and be defined as the union of all translates of $\rho\mathbf{C}$ that are covered by $\mathbf{K}$.

Now, take a fixed integer $m\ge 1$. On the one hand, if $\rho>0$ is sufficiently small, then $ w_{\mathbf{C}}({\mathbf{K}}^{\rho\mathbf{C}})>m\rho$. On the other hand, $ w_{\mathbf{C}}({\mathbf{K}}^{r\mathbf{C}})=r\le mr$. As $ w_{\mathbf{C}}({\mathbf{K}}^{\rho\mathbf{C}})$ is a decreasing continuous function of $\rho>0$ and $m\rho$ is a strictly increasing continuous function of $\rho$, there exists a uniquely determined $\rho>0$ such that
$$ w_{\mathbf{C}}({\mathbf{K}}^{\rho\mathbf{C}})=m\rho.$$ This uniquely determined $\rho$ is called the {\it $m$th successive $\mathbf{C}$-inradius of $\mathbf{K}$} and is denoted by $r_{\mathbf{C}}(\mathbf{K} , m)$.

Now, the two equivalent versions of Conjecture \ref{Bang-conjecture} and Conjecture \ref{Alexander--Davenport} introduced in \cite{BeBe96} can be phrased as follows.

\begin{conjecture}\label{Bezdek--Bezdek-1}
If a convex body $\mathbf{K}$ in $\mathbb{E}^d, d\ge 2$ is covered by the planks $\mathbf{P}_1, \mathbf{P}_2,$ $\dots , \mathbf{P}_n$, then $\sum_{i=1}^n w_{\mathbf{C}}(\mathbf{P}_i)\ge
w_{\mathbf{C}}(\mathbf{K})$ for any convex body $\mathbf{C}$ in $\mathbb{E}^d$.
\end{conjecture}

\begin{conjecture}\label{Bezdek--Bezdek-2}
Let $\mathbf{K}$ and $\mathbf{C}$ be convex bodies in $\mathbb{E}^d, d\ge 2$. If $\mathbf{K}$ is sliced by $n-1$ hyperplanes, then the minimum of the greatest $\mathbf{C}$-inradius of the pieces is equal to the $n$th successive $\mathbf{C}$-inradius of $\mathbf{K}$, i.e., it is $r_{\mathbf{C}}(\mathbf{K} , n)$.
\end{conjecture}

Recall that Theorem \ref{Ball-plank-th} gives a proof of (Conjecture \ref{Bezdek--Bezdek-1} as well as) Conjecture \ref{Bezdek--Bezdek-2} for centrally symmetric convex bodies $\mathbf{K}$ in $\mathbb{E}^d, d\ge 2$ (with $\mathbf{C}$ being an arbitrary convex body in $\mathbb{E}^d, d\ge 2$). Another approach that leads to a partial solution of Conjecture \ref{Bezdek--Bezdek-2} was published in \cite{BeBe96}. Namely, in that paper A. Bezdek and the author proved the following theorem that (under the condition that $\mathbf{C}$ is a ball) answers a question raised by Conway (\cite{BeBe95}) as well as proves Conjecture \ref{Bezdek--Bezdek-2} for successive hyperplane cuts.

\begin{theorem}\label{Bezdek-Bezdek-Conway}
Let $\mathbf{K}$ and $\mathbf{C}$ be convex bodies in $\mathbb{E}^d$, $d\ge 2$. If $\mathbf{K}$ is sliced into $n\ge 1$ pieces by $n-1$ successive hyperplane cuts (i.e., when each cut divides one piece), then the minimum of the greatest $\mathbf{C}$-inradius of the pieces is the $n$th successive $\mathbf{C}$-inradius of $\mathbf{K}$ (i.e., $r_{\mathbf{C}}(\mathbf{K} , n)$). An optimal partition is achieved by $n-1$ parallel hyperplane cuts equally spaced along the minimal $\mathbf{C}$-width of the $r_{\mathbf{C}}(\mathbf{K} , n)\mathbf{C}$-rounded body of $\mathbf{K}$.
\end{theorem}

Akopyan and Karasev (\cite{AkKa12}) just very recently have proved a related partial result on Conjecture~\ref{Bezdek--Bezdek-1}. Their theorem is based on a nice generalization of successive hyperplane cuts. The more exact details are as follows. Under the {\it convex partition} $\mathbf{V}_1\cup\mathbf{V}_2\cup\dots\cup\mathbf{V}_n$ of $\mathbb{E}^d$ we understand the family $\mathbf{V}_1, \mathbf{V}_2, \dots, \mathbf{V}_n$ of closed convex sets having pairwise disjoint non-empty interiors in $\mathbb{E}^d$ with $\mathbf{V}_1\cup\mathbf{V}_2\cup\dots\cup\mathbf{V}_n=\mathbb{E}^d$. Then we say that the convex partition $\mathbf{V}_1\cup\mathbf{V}_2\cup\dots\cup\mathbf{V}_n$ of $\mathbb{E}^d$
is an {\it inductive partition} of $\mathbb{E}^d$ if for any $1\le i\le n$, there exists an inductive partition $\mathbf{W}_1\cup\dots\cup\mathbf{W}_{i-1}\cup\mathbf{W}_{i+1}\cup\dots\cup\mathbf{W}_n$ of $\mathbb{E}^d$ such that $\mathbf{V}_j\subset\mathbf{W}_j$ for all $j\neq i$. A partition into one part $\mathbf{V}_1=\mathbb{E}^d$ is assumed to be inductive. We note that if $\mathbb{E}^d$ is sliced into $n$ pieces by $n-1$ successive hyperplane cuts (i.e., when each cut divides one piece), then the pieces generate an inductive partition of $\mathbb{E}^d$. Also, the Voronoi cells of finitely many points of $\mathbb{E}^d$ generate an inductive partition of $\mathbb{E}^d$. Now, the main theorem of \cite{AkKa12} can be phrased as follows.

\begin{theorem}\label{Akopyan-Karasev}
Let $\mathbf{K}$ and $\mathbf{C}$ be convex bodies in $\mathbb{E}^d, d\ge 2$ and let $\mathbf{V}_1\cup\mathbf{V}_2\cup\dots\cup\mathbf{V}_n$ be an inductive partition of $\mathbb{E}^d$ such that ${\rm int}(\mathbf{V}_i\cap\mathbf{K})\neq\emptyset$ for all $1\le i\le n$. Then 
$\sum_{i=1}^{n}r_{\mathbf{C}}(\mathbf{V}_i\cap\mathbf{K}, 1)\ge r_{\mathbf{C}}(\mathbf{K} , 1)$. 
\end{theorem} 

\section{Extensions to Successive Inradii}

First, we state the following stronger version of Theorem~\ref{Bezdek-Bezdek-Conway}. Its proof is an extension of the proof of Theorem~\ref{Bezdek-Bezdek-Conway} published in \cite{BeBe96}. 

\begin{theorem}\label{Bezdek-Bezdek-Conway-generalized}
Let $\mathbf{K}$ and $\mathbf{C}$ be convex bodies in $\mathbb{E}^d$, $d\ge 2$ and let $m$ be a positive integer. If $\mathbf{K}$ is sliced into $n\ge 1$ pieces by $n-1$ successive hyperplane cuts (i.e., when each cut divides one piece), then the minimum of the greatest $m$th successive $\mathbf{C}$-inradius of the pieces is the $(mn)$th successive $\mathbf{C}$-inradius of $\mathbf{K}$ (i.e., $r_{\mathbf{C}}(\mathbf{K} , mn)$). An optimal partition is achieved by $n-1$ parallel hyperplane cuts equally spaced along the minimal $\mathbf{C}$-width of the $r_{\mathbf{C}}(\mathbf{K} , mn)\mathbf{C}$-rounded body of $\mathbf{K}$.
\end{theorem}

Second, the method of Akopyan and Karasev (\cite{AkKa12}) can be extended to prove the following stronger version of Theorem \ref{Akopyan-Karasev}. In fact, that approach extends also the relavant additional theorems of Akopyan and Karasev stated in \cite{AkKa12} and used in their proof of Theorem~\ref{Akopyan-Karasev}. However, in this paper following the recommendation of the referee, we derive the next theorem directly from Theorem~\ref{Akopyan-Karasev}. 

\begin{theorem}\label{Akopyan-Karasev-Bezdek}
Let $\mathbf{K}$ and $\mathbf{C}$ be convex bodies in $\mathbb{E}^d, d\ge 2$  and let $m$ be a positive integer. If $\mathbf{V}_1\cup\mathbf{V}_2\cup\dots\cup\mathbf{V}_n$ is an inductive partition of $\mathbb{E}^d$ such that ${\rm int}(\mathbf{V}_i\cap\mathbf{K})\neq\emptyset$ for all $1\le i\le n$, then 
$\sum_{i=1}^{n}r_{\mathbf{C}}(\mathbf{V}_i\cap\mathbf{K}, m)\ge r_{\mathbf{C}}(\mathbf{K} , m)$.
\end{theorem}

\begin{corollary}\label{corollary of Akopyan-Karasev-Bezdek}
Let $\mathbf{K}$ and $\mathbf{C}$ be convex bodies in $\mathbb{E}^d, d\ge 2$. If $\mathbf{V}_1\cup\mathbf{V}_2\cup\dots\cup\mathbf{V}_n$ is an inductive partition of $\mathbb{E}^d$ such that ${\rm int}(\mathbf{V}_i\cap\mathbf{K})\neq\emptyset$ for all $1\le i\le n$, then 
$\sum_{i=1}^{n}w_{\mathbf{C}}(\mathbf{V}_i\cap\mathbf{K})\ge w_{\mathbf{C}}(\mathbf{K})$.
\end{corollary}

For the sake of completeness we mention that in two dimensions one can state a bit more. Namely, recall that Akopyan and Karasev (\cite{AkKa12}) proved the following: Let $\mathbf{K}$ and $\mathbf{C}$ be convex bodies in $\mathbb{E}^2$ and let $\mathbf{V}_1\cup\mathbf{V}_2\cup\dots\cup\mathbf{V}_n=\mathbf{K}$ be a partition of $\mathbf{K}$ into convex bodies $\mathbf{V}_i$, $1\le i\le n$. Then $\sum_{i=1}^{n}r_{\mathbf{C}}(\mathbf{V}_i, 1)\ge r_{\mathbf{C}}(\mathbf{K} , 1)$. Now, exactly the same way as Theorem~\ref{Akopyan-Karasev-Bezdek} is derived from Theorem~\ref{Akopyan-Karasev}, it follows that $\sum_{i=1}^{n}r_{\mathbf{C}}(\mathbf{V}_i, m)\ge r_{\mathbf{C}}(\mathbf{K} , m)$ holds for any positive integer $m$.

Finally, we close this section stating that Conjectures~\ref{Bang-conjecture}, ~\ref{Alexander--Davenport}, ~\ref{Bezdek--Bezdek-1}, and ~\ref{Bezdek--Bezdek-2} are all equivalent to the following two conjectures:

\begin{conjecture}\label{Bezdek--Bezdek-11}
Let $\mathbf{K}$ and $\mathbf{C}$ be convex bodies in $\mathbb{E}^d, d\ge 2$  and let $m$ be a positive integer. If $\mathbf{K}$  is covered by the planks $\mathbf{P}_1, \mathbf{P}_2,$ $\dots , \mathbf{P}_n$ in $\mathbb{E}^d$, then $\sum_{i=1}^n r_{\mathbf{C}}(\mathbf{P}_i , m)\ge r_{\mathbf{C}}(\mathbf{K} , m)$ or equivalently, $\sum_{i=1}^n w_{\mathbf{C}}(\mathbf{P}_i)\ge m r_{\mathbf{C}}(\mathbf{K} , m)$.
\end{conjecture}

\begin{conjecture}\label{Bezdek--Bezdek-22} 
Let $\mathbf{K}$ and $\mathbf{C}$ be convex bodies in $\mathbb{E}^d, d\ge 2$  and let the positive integer $m$ be given. If $\mathbf{K}$ is sliced by $n-1$ hyperplanes, then the minimum of the greatest $m$th successive $\mathbf{C}$-inradius of the pieces is the $(mn)$th successive $\mathbf{C}$-inradius of $\mathbf{K}$, i.e., it is $r_{\mathbf{C}}(\mathbf{K} , mn)$.
\end{conjecture}

In the rest of the paper we prove the claims of this section.

\section{Proof of Theorem \ref{Bezdek-Bezdek-Conway-generalized}}

\subsection{On Coverings of Convex Bodies by Two Planks}

On the one hand, the following statement is an extension to higher dimensions of Theorem 4 in \cite{Al68}. On the other hand, the proof presented below is based on Theorem 4 of \cite{Al68}.

\begin{lemma}\label{Bezdek-Bezdek-Alexander-inequality}
If a convex body $\mathbf{K}$ in $\mathbb{E}^d, d\ge 2$ is covered by the planks $\mathbf{P}_1$ and $\mathbf{P}_2$, then $w_{\mathbf{C}}(\mathbf{P}_1)+w_{\mathbf{C}}(\mathbf{P}_2)\ge w_{\mathbf{C}}(\mathbf{K})$ for any convex body $\mathbf{C}$ in $\mathbb{E}^d$.
\end{lemma}

\proof
Let $H_1$ (resp., $H_2$) be one of the two hyperplanes which bound the plank $\mathbf{P}_1$ (resp., $\mathbf{P}_2$). If $H_1$ and $H_2$ are translates of each other, then the claim is obviously true. Thus, without loss of generality we may assume that $L:=H_1\cap H_2$ is a $(d-2)$-dimensional affine subspace of $\mathbb{E}^d$. Let $\mathbb{E}^2$ be the $2$-dimensional linear subspace of $\mathbb{E}^d$ that is orthogonal to $L$. If $(\cdot)'$ denotes the (orthogonal) projection of $\mathbb{E}^d$ parallel to $L$ onto $\mathbb{E}^2$, then obviously, $w_{\mathbf{C'}}(\mathbf{P}_1')=w_{\mathbf{C}}(\mathbf{P}_1)$, $w_{\mathbf{C'}}(\mathbf{P}_2')=w_{\mathbf{C}}(\mathbf{P}_2)$ and $w_{\mathbf{C'}}(\mathbf{K'})\ge w_{\mathbf{C}}(\mathbf{K})$. Thus, it is sufficient to prove that
$$
w_{\mathbf{C'}}(\mathbf{P}_1')+w_{\mathbf{C'}}(\mathbf{P}_2')\ge w_{\mathbf{C'}}(\mathbf{K'}).
$$
In other words, it is sufficient to prove Lemma \ref{Bezdek-Bezdek-Alexander-inequality} for $d=2$. Hence, in the rest of the proof, $\mathbf{K}, \mathbf{C}, \mathbf{P}_1, \mathbf{P}_2, H_1 $, and $H_2$ mean the sets introduced and defined above, however, for $d=2$. Now, we can make the following easy observation
$$
w_{\mathbf{C}}(\mathbf{P}_1)+w_{\mathbf{C}}(\mathbf{P}_2)=\frac{w(\mathbf{P}_1)}{w(\mathbf{C}, H_1)}+\frac{w(\mathbf{P}_2)}{w(\mathbf{C}, H_2)}
$$
$$
=\frac{w(\mathbf{P}_1)}{w(\mathbf{K}, H_1)}\frac{w(\mathbf{K}, H_1)}{w(\mathbf{C}, H_1)}+\frac{w(\mathbf{P}_2)}{w(\mathbf{K}, H_2)}\frac{w(\mathbf{K}, H_2)}{w(\mathbf{C}, H_2)}
$$
$$
\ge \left(\frac{w(\mathbf{P}_1)}{w(\mathbf{K}, H_1)}+ \frac{w(\mathbf{P}_2)}{w(\mathbf{K}, H_2)} \right)w_{\mathbf{C}}(\mathbf{K})
$$
$$
=\left(w_{\mathbf{K}}(\mathbf{P}_1)+ w_{\mathbf{K}}(\mathbf{P}_2)\right)w_{\mathbf{C}}(\mathbf{K}).
$$
Then recall that Theorem 4 in \cite{Al68} states that if a convex set in the plane is covered by two planks, then the sum of their relative widths is at least $1$. Thus, using our terminology, we have that $w_{\mathbf{K}}(\mathbf{P}_1)+ w_{\mathbf{K}}(\mathbf{P}_2)\ge 1$, finishing the proof of Lemma \ref{Bezdek-Bezdek-Alexander-inequality}.
\endproof

\subsection{Minimizing the Greatest $m$th Successive $\mathbf{C}$-Inradius}

Let $\mathbf{K}$ and $\mathbf{C}$ be convex bodies in $\mathbb{E}^d$, $d\ge 2$. We prove Theorem \ref{Bezdek-Bezdek-Conway-generalized} by induction on $n$. It is trivial to check the claim for $n=1$. So, let $n\ge 2$ be given and assume that Theorem \ref{Bezdek-Bezdek-Conway-generalized} holds for at most $n-2$  successive hyperplane cuts and based on that we show that it holds for $n-1$ successive hyperplane cuts as well. The details are as follows.

Let $H_1, \dots, H_{n-1}$ denote the hyperplanes of the $n-1$ successive hyperplane cuts that slice $\mathbf{K}$ into $n$ pieces such that the greatest $m$th successive $\mathbf{C}$-inradius of the pieces is the smallest possible say, $\rho$. Then take the first cut $H_1$ that slices $\mathbf{K}$ into the pieces $\mathbf{K}_1$ and $\mathbf{K}_2$ such that $\mathbf{K}_1$ (resp., $\mathbf{K}_2$) is sliced into $n_1$ (resp., $n_2$) pieces by the successive hyperplane cuts $H_2, \dots, H_{n-1}$, where $n=n_1+n_2$. The induction hypothesis implies that $\rho\ge r_{\mathbf{C}}(\mathbf{K}_1 , mn_1)=:\rho_1$ and $\rho\ge r_{\mathbf{C}}(\mathbf{K}_2 , mn_2)=:\rho_2$ and therefore
\begin{equation}\label{induction-hypothesis-1}
w_{\mathbf{C}}({\mathbf{K}_1}^{\rho\mathbf{C}})\le w_{\mathbf{C}}({\mathbf{K}_1}^{\rho_1\mathbf{C}})=mn_1\rho_1\le mn_1\rho ;
\end{equation}
moreover,
\begin{equation}\label{induction-hypothesis-2}
w_{\mathbf{C}}({\mathbf{K}_2}^{\rho\mathbf{C}})\le w_{\mathbf{C}}({\mathbf{K}_2}^{\rho_2\mathbf{C}})=mn_2\rho_2\le mn_2\rho .
\end{equation}
Now, we need to define the following set. 

\begin{definition} Assume that the origin $\mathbf{o}$ of $\mathbb{E}^d$ belongs to the interior of the convex body $\mathbf{C}\subset\mathbb{E}^d$. Consider all translates of $\rho\mathbf{C}$ which are contained in the convex body $\mathbf{K}\subset\mathbb{E}^d$. The set of points in the translates of $\rho\mathbf{C}$ that correspond to $\mathbf{o}$ form a convex set called the inner $\rho\mathbf{C}$-parallel body of $\mathbf{K}$ denoted by $\mathbf{K}_{-\rho\mathbf{C}}$. 
\end{definition}

Clearly,
$$
(\mathbf{K}_1)_{-\rho\mathbf{C}}\cup(\mathbf{K}_2)_{-\rho\mathbf{C}}\subset\mathbf{K}_{-\rho\mathbf{C}} \ {\rm with}\ (\mathbf{K}_1)_{-\rho\mathbf{C}}\cap(\mathbf{K}_2)_{-\rho\mathbf{C}}=\emptyset .
$$
Also, it is easy to see that there is a plank $\mathbf{P}$ with $w_{\mathbf{C}}(\mathbf{P})=\rho$ such that it is parallel to $H_1$ and contains $H_1$ in its interior; moreover,
$$
\mathbf{K}_{-\rho\mathbf{C}}\subset(\mathbf{K}_1)_{-\rho\mathbf{C}}\cup(\mathbf{K}_2)_{-\rho\mathbf{C}}\cup \mathbf{P} .
$$

\noindent Now, let $H_1^+$ (resp., $H_1^-$) be the closed halfspace of $\mathbb{E}^d$ bounded by $H_1$ and containing $\mathbf{K}_1$ (resp., $\mathbf{K}_2$) and let
$\mathbf{P}^+:=\mathbf{P}\cap H_1^+$ (resp., $\mathbf{P}^-:=\mathbf{P}\cap H_1^-$). Moreover, let $\mathbf{K}_{-\rho\mathbf{C}}^+:=\mathbf{K}_{-\rho\mathbf{C}}\cap H_1^+$
(resp., $\mathbf{K}_{-\rho\mathbf{C}}^-:=\mathbf{K}_{-\rho\mathbf{C}}\cap H_1^-$). Hence, applying Lemma \ref{Bezdek-Bezdek-Alexander-inequality} to $\mathbf{K}_{-\rho\mathbf{C}}$ partitioned into $\mathbf{K}_{-\rho\mathbf{C}}^+\cup\mathbf{K}_{-\rho\mathbf{C}}^-$ and to $\mathbf{K}_{-\rho\mathbf{C}}^+$ covered by the plank $\mathbf{P}^+$ and the plank generated by the minimal $\mathbf{C}$-width of $(\mathbf{K}_1)_{-\rho\mathbf{C}}$ as well as to $\mathbf{K}_{-\rho\mathbf{C}}^-$ covered by the plank $\mathbf{P}^-$  and the plank generated by the minimal $\mathbf{C}$-width of $(\mathbf{K}_2)_{-\rho\mathbf{C}}$ we get that 
\begin{equation}\label{width-inequality-for-inner-parallel-bodies}
w_{\mathbf{C}}\left(\mathbf{K}_{-\rho\mathbf{C}}\right)\le w_{\mathbf{C}}\left(\mathbf{K}_{-\rho\mathbf{C}}^+\right)+w_{\mathbf{C}}\left(\mathbf{K}_{-\rho\mathbf{C}}^-\right) 
\le w_{\mathbf{C}}\left((\mathbf{K}_1)_{-\rho\mathbf{C}}\right)+\rho+w_{\mathbf{C}}\left((\mathbf{K}_2)_{-\rho\mathbf{C}}\right).
\end{equation}

\noindent By definition $w_{\mathbf{C}}\left((\mathbf{K}_1)_{-\rho\mathbf{C}}\right)$ $=w_{\mathbf{C}}({\mathbf{K}_1}^{\rho\mathbf{C}})-\rho$, $w_{\mathbf{C}}\left((\mathbf{K}_2)_{-\rho\mathbf{C}}\right)$ $=w_{\mathbf{C}}({\mathbf{K}_2}^{\rho\mathbf{C}})-\rho$ and $w_{\mathbf{C}}\left(\mathbf{K}_{-\rho\mathbf{C}}\right)=w_{\mathbf{C}}(\mathbf{K}^{\rho\mathbf{C}})-\rho$. Hence, (\ref{width-inequality-for-inner-parallel-bodies}) is equivalent to
\begin{equation}\label{width-inequality-for-rounded-bodies}
w_{\mathbf{C}}(\mathbf{K}^{\rho\mathbf{C}})\le w_{\mathbf{C}}({\mathbf{K}_1}^{\rho\mathbf{C}})+w_{\mathbf{C}}({\mathbf{K}_2}^{\rho\mathbf{C}}).
\end{equation} 
Finally, (\ref{induction-hypothesis-1}),(\ref{induction-hypothesis-2}), and (\ref{width-inequality-for-rounded-bodies}) yield that
\begin{equation}\label{final-inductive-inequality}
w_{\mathbf{C}}(\mathbf{K}^{\rho\mathbf{C}})\le mn_1\rho+mn_2\rho=mn\rho.
\end{equation}
Thus, (\ref{final-inductive-inequality}) clearly implies that $r_{\mathbf{C}}(\mathbf{K} , mn)\le \rho$. As the case, when the optimal partition is achieved, follows directly from the definition of the $mn$th successive $\mathbf{C}$-inradius of $\mathbf{K}$, the proof of Theorem \ref{Bezdek-Bezdek-Conway-generalized} is complete.

\section{Proof of Theorem~\ref{Akopyan-Karasev-Bezdek}}

Let $\mathbf{K}$ and $\mathbf{C}$ be convex bodies in $\mathbb{E}^d, d\ge 2$  and let $m$ be a positive integer. It follows from the definition of $r_{\mathbf{C}}(\mathbf{K} , m)$ that $r_{\mathbf{C}}(\mathbf{K} , m)$ is a translation invariant, positively $1$-homogeneous, inclusion-monotone functional over the family of convex bodies $\mathbf{K}$ in $\mathbb{E}^d$ for any fixed $\mathbf{C} $ and $m$. On the other hand, if $\mathbf{V}_1\cup\mathbf{V}_2\cup\dots\cup\mathbf{V}_n$ is an inductive partition of $\mathbb{E}^d$ such that ${\rm int}(\mathbf{V}_i\cap\mathbf{K})\neq\emptyset$ for all $1\le i\le n$, then 
Theorem~\ref{Akopyan-Karasev} applied to $\mathbf{C}=\mathbf{K}$ yields the existence of translation vectors $\mathbf{t}_1, \mathbf{t}_2, \dots , \mathbf{t}_n$ and positive reals $\mu_1,\mu_2,\dots ,\mu_n$ such that $\mathbf{t}_i+\mu_i \mathbf{K}\subset \mathbf{V}_i\cap\mathbf{K}$ for all $1\le i\le n$ satisfying $\sum_{i=1}^{n}\mu_i\ge 1$. Therefore $$r_{\mathbf{C}}(\mathbf{V}_i\cap\mathbf{K} , m)\ge r_{\mathbf{C}}(\mathbf{t}_i+\mu_i \mathbf{K}  , m)=r_{\mathbf{C}}(\mu_i \mathbf{K}  , m)=\mu_i r_{\mathbf{C}}(\mathbf{K} , m)$$ holds for all $1\le i\le n$, finishing the proof of Theorem~\ref{Akopyan-Karasev-Bezdek}.

\section{Proof of Corollary~\ref{corollary of Akopyan-Karasev-Bezdek}}

Let $1\le m_1\le m_2$ be positive integers. Recall that if $\rho_1$ (resp., $\rho_2$) denotes the $m_1$th (resp., $m_2$th) successive $\mathbf{C}$-inradius of $\mathbf{K}$, then by definition 
$ w_{\mathbf{C}}({\mathbf{K}}^{\rho_1\mathbf{C}})=m_1\rho_1$ (resp., $ w_{\mathbf{C}}({\mathbf{K}}^{\rho_2\mathbf{C}})=m_2\rho_2$). As $ w_{\mathbf{C}}({\mathbf{K}}^{\rho\mathbf{C}})$ is a decreasing continuous function of $\rho>0$, it follows that
$$m_1r_{\mathbf{C}}(\mathbf{K} , m_1)=m_1\rho_1\le m_2\rho_2=m_2r_{\mathbf{C}}(\mathbf{K} , m_2)\ .$$
Thus,  the sequence $mr_{\mathbf{C}}(\mathbf{K} , m), m=1,2,\dots$ is an increasing one with  
$$\lim_{m\to+\infty}mr_{\mathbf{C}}(\mathbf{K} , m)=w_{\mathbf{C}}(\mathbf{K})\ .$$
Hence, Corollary~\ref{corollary of Akopyan-Karasev-Bezdek} follows from Theorem \ref{Akopyan-Karasev-Bezdek}.

\section{The equivalence of Conjectures~\ref{Bang-conjecture}, ~\ref{Alexander--Davenport}, ~\ref{Bezdek--Bezdek-1}, ~\ref{Bezdek--Bezdek-2}, ~\ref{Bezdek--Bezdek-11}, and ~\ref{Bezdek--Bezdek-22}}

Recall that according to \cite{BeBe96} Conjectures~\ref{Bang-conjecture}, ~\ref{Alexander--Davenport}, ~\ref{Bezdek--Bezdek-1}, and ~\ref{Bezdek--Bezdek-2} are equivalent
to each other. So, it is sufficent to show that Conjecture~\ref{Bezdek--Bezdek-1} implies Conjecture~\ref{Bezdek--Bezdek-11} and Conjecture~\ref{Bezdek--Bezdek-11} implies Conjecture~\ref{Bezdek--Bezdek-22} moreover,  Conjecture~\ref{Bezdek--Bezdek-22} implies Conjecture~\ref{Bezdek--Bezdek-2}.

As according to the previous section the sequence $mr_{\mathbf{C}}(\mathbf{K} , m), m=1,2,\dots$ is an increasing one with  
$\lim_{m\to+\infty}mr_{\mathbf{C}}(\mathbf{K} , m)=w_{\mathbf{C}}(\mathbf{K})$ therefore Conjecture~\ref{Bezdek--Bezdek-1} implies Conjecture~\ref{Bezdek--Bezdek-11}. Next, it is obvious that Conjecture~\ref{Bezdek--Bezdek-22} implies Conjecture~\ref{Bezdek--Bezdek-2}. So, we are left to show that Conjecture~\ref{Bezdek--Bezdek-11} implies Conjecture~\ref{Bezdek--Bezdek-22}.  In order to do so we introduce the following equivalent description for $r_{\mathbf{C}}(\mathbf{K} , m)$. If $\mathbf{C}$ is a convex body in $\mathbb{E}^d$, then $$\mathbf{t}+\mathbf{C}, \mathbf{t}+\lambda_2\mathbf{v}+\mathbf{C}, \dots , \mathbf{t}+\lambda_m\mathbf{v}+\mathbf{C}$$ is called a {\it linear packing} of $m$ translates of $\mathbf{C}$ positioned parallel to the 
line $\{\lambda\mathbf{v}\ |\ \lambda\in\mathbb{R}\}$ with direction vector $\mathbf{v}\neq\mathbf{o}$ if the $m$ translates of $\mathbf{C}$ are pairwise non-overlapping, i.e., if
$$( \mathbf{t}+\lambda_i\mathbf{v}+{\rm int}\mathbf{C}) \cap ( \mathbf{t}+\lambda_j\mathbf{v}+{\rm int}\mathbf{C})=\emptyset$$ holds for all $1\le i\neq j\le m$ (with $\lambda_1=0$). 
Furthermore, the line $l\subset \mathbb{E}^d$ passing through the origin $\mathbf{o}$ of $\mathbb{E}^d$ is called a {\it separating direction} for the linear packing 
$$\mathbf{t}+\mathbf{C}, \mathbf{t}+\lambda_2\mathbf{v}+\mathbf{C}, \dots , \mathbf{t}+\lambda_m\mathbf{v}+\mathbf{C}$$
if 
$${\rm Pr}_{l}( \mathbf{t}+\mathbf{C}), {\rm Pr}_{l}( \mathbf{t}+\lambda_2\mathbf{v}+\mathbf{C}), \dots , {\rm Pr}_{l}( \mathbf{t}+\lambda_m\mathbf{v}+\mathbf{C})$$ 
are pairwise non-overlapping intervals on $l$, where ${\rm Pr}_l: \mathbb{E}^d\to l$ denotes the orthogonal projection of $\mathbb{E}^d$ onto $l$. It is easy to see that every linear packing 
$$\mathbf{t}+\mathbf{C}, \mathbf{t}+\lambda_2\mathbf{v}+\mathbf{C}, \dots , \mathbf{t}+\lambda_m\mathbf{v}+\mathbf{C}$$
possesses at least one separating direction in $\mathbb{E}^d$. Finally, let $\mathbf{K}$ be a convex body in $\mathbb{E}^d$ and let $m\ge 1$ be a positive integer. Then let $\overline{\rho}>0$ 
be the largest positive real with the following property: for every line $l$ passing through the origin $\mathbf{o}$ in $\mathbb{E}^d$ there exists a linear packing of $m$ translates of $\overline{\rho}\mathbf{C}$ lying in $\mathbf{K}$ and having $l$ as a separating direction.
It is straightforward to show that $$\overline{\rho}=r_{\mathbf{C}}(\mathbf{K} , m).$$ 
Now, let $\mathbf{K}$ and $\mathbf{C}$ be convex bodies in $\mathbb{E}^d, d\ge 2$  and let the positive integer $m$ be given. Assume that the origin $\mathbf{o}$ of $\mathbb{E}^d$ lies in the interior of $\mathbf{C}$. Furthermore, assume that $\mathbf{K}$ is sliced by $n-1$ hyperplanes say, $H_1, H_2, \dots , H_{n-1}$ and let $\rho$ be the greatest $m$th successive $\mathbf{C}$-inradius of the pieces of $\mathbf{K}$ obtained in this way. Then let $\mathbf{P}_i:=\bigcup_{\mathbf{p}\in H_i}\left(\mathbf{p}+(-m\rho)\mathbf{C}\right)$, $1\le i\le n-1$. Based on the above description of $m$th successive $\mathbf{C}$-inradii, it is easy to see that $\mathbf{K}_{-m\rho\mathbf{C}}\subset    \bigcup_{i=1}^{n-1}\mathbf{P}_i$ with $w_{\mathbf{C}}(\mathbf{P}_i )=m\rho$ for all $1\le i\le n-1$. Thus, Conjecture~\ref{Bezdek--Bezdek-11} implies that
$(n-1)m\rho=\sum_{i=1}^{n-1} w_{\mathbf{C}}(\mathbf{P}_i )\ge m r_{\mathbf{C}}(  \mathbf{K}_{-m\rho\mathbf{C}}  , m)=m\left(r_{\mathbf{C}}(\mathbf{K}^{\rho\mathbf{C}}, m)-\rho\right)$ and so, $mn\rho\ge w_{\mathbf{C}}(\mathbf{K}^{\rho\mathbf{C}})$. Hence,
$\rho\ge r_{\mathbf{C}}(\mathbf{K} , mn)$ finishing the proof of Conjecture~\ref{Bezdek--Bezdek-22}.

\section{Conclusion}

Theorems~\ref{Akopyan-Karasev} and ~\ref{Akopyan-Karasev-Bezdek} have covering analogues. Namely recall that Akopyan and Karasev (\cite{AkKa12}) introduced the following definition. Under the {\it convex covering} $\mathbf{V}_1\cup\mathbf{V}_2\cup\dots\cup\mathbf{V}_n$ of $\mathbb{E}^d$ we understand the family $\mathbf{V}_1, \mathbf{V}_2, \dots, \mathbf{V}_n$ of closed convex sets in $\mathbb{E}^d$ with $\mathbf{V}_1\cup\mathbf{V}_2\cup\dots\cup\mathbf{V}_n=\mathbb{E}^d$. Then we say that the convex covering $\mathbf{V}_1\cup\mathbf{V}_2\cup\dots\cup\mathbf{V}_n$ of $\mathbb{E}^d$
is an {\it inductive covering} of $\mathbb{E}^d$ if for any $1\le i\le n$, there exists an inductive covering $\mathbf{W}_1\cup\dots\cup\mathbf{W}_{i-1}\cup\mathbf{W}_{i+1}\cup\dots\cup\mathbf{W}_n$ of $\mathbb{E}^d$ such that $\mathbf{W}_j\subset\mathbf{V}_j\cup \mathbf{V}_i$ for all $j\neq i$. A covering by one set $\mathbf{V}_1=\mathbb{E}^d$ is assumed to be inductive. \cite{AkKa12} proves that if $\mathbf{K}$ and $\mathbf{C}$ are convex bodies in $\mathbb{E}^d, d\ge 2$ and $\mathbf{V}_1\cup\mathbf{V}_2\cup\dots\cup\mathbf{V}_n$ is an inductive covering of $\mathbb{E}^d$ such that ${\rm int}(\mathbf{V}_i\cap\mathbf{K})\neq\emptyset$ for all $1\le i\le n$, then 
$\sum_{i=1}^{n}r_{\mathbf{C}}(\mathbf{V}_i\cap\mathbf{K}, 1)\ge r_{\mathbf{C}}(\mathbf{K} , 1)$. Now, exactly the same way as Theorem~\ref{Akopyan-Karasev-Bezdek} is derived from Theorem~\ref{Akopyan-Karasev}, it follows that
\begin{equation} \label{Akopyan-Karasev-II}
\sum_{i=1}^{n}r_{\mathbf{C}}(\mathbf{V}_i\cap\mathbf{K}, m)\ge r_{\mathbf{C}}(\mathbf{K} , m) 
\end{equation}
holds for any positive integer $m$. This raises the following rather natural question (see also Conjecture~\ref{Bezdek--Bezdek-11}).

\begin{problem}
Let $\mathbf{K}$ and $\mathbf{C}$ be convex bodies in $\mathbb E^d$, $d\geq 2$ and let $m$ be a positive integer.
Prove or disprove that if $\mathbf{V}_1 \cup \mathbf{V}_2 \cup \ldots \cup \mathbf{V}_n$ is
a convex partition (resp., covering) of $\mathbb E^d$ such that ${\rm int}(\mathbf{V}_i\cap\mathbf{K})\neq\emptyset$ for all $1\le i\le n$, then
$\sum_{i=1}^{n}r_{\mathbf{C}}(\mathbf{V}_i\cap\mathbf{K}, m)\ge r_{\mathbf{C}}(\mathbf{K} , m)$.
\end{problem}

Next observe that (\ref{Akopyan-Karasev-II}) implies in a straightforward way that if $\mathbf{K}$ and $\mathbf{C}$ are convex bodies in $\mathbb E^d$ and $\mathbf{V}_1 \cup \mathbf{V}_2 \cup \ldots \cup \mathbf{V}_n$ is an inductive covering of $\mathbb E^d$ such that ${\rm int}(\mathbf{V}_i\cap\mathbf{K})\neq\emptyset$ for all $1\le i\le n$, then the greatest
$m$th successive $\mathbf{C}$-inradius of the pieces $\mathbf{V}_i\cap\mathbf{K}, i=1, 2, \dots , n$ is at least $\frac{1}{n}r_{\mathbf{C}}(\mathbf{K} , m)$. As the sequence $mr_{\mathbf{C}}(\mathbf{K} , m), m=1,2,\dots$ is an increasing one, therefore $\frac{1}{n}r_{\mathbf{C}}(\mathbf{K} , m)\le r_{\mathbf{C}}(\mathbf{K} , mn)$ raising the following question (see also Conjecture~\ref{Bezdek--Bezdek-22}). 

\begin{problem}
Let $\mathbf{K}$ and $\mathbf{C}$ be convex bodies in $\mathbb E^d$, $d\geq 2$ and let $m$ be a positive integer.
Prove or disprove that if $\mathbf{V}_1 \cup \mathbf{V}_2 \cup \ldots \cup \mathbf{V}_n$ is
a convex partition (resp., covering) of $\mathbb E^d$ such that ${\rm int}(\mathbf{V}_i\cap\mathbf{K})\neq\emptyset$ for all $1\le i\le n$, then
 the greatest $m$th successive $\mathbf{C}$-inradius of the pieces $\mathbf{V}_i\cap\mathbf{K}, i=1, 2, \dots , n$ is at least
$r_{\mathbf{C}}(\mathbf{K} , mn)$.
\end{problem}

\bibliographystyle{amsplain}

\end{document}